\RequirePackage{ifpdf}
\ifpdf 
\documentclass[pdftex]{sigma}
\else
\documentclass{sigma}
\fi

\usepackage[mathscr]{eucal}

\usepackage[all]{xy}
\CompileMatrices

\numberwithin{equation}{section}

\newfont{\germ}{eufm10}

\newcommand\et[1]{\tilde{e}_{#1}}
\newcommand\ft[1]{\tilde{f}_{#1}}
\newcommand\geh{\goth{g}}

\newcommand\goth[1]{\mbox{\germ #1}}

\newcommand\La{\Lambda}
\newcommand\la{\lambda}
\newcommand\lw[1]{\lower.4mm\hbox{${#1}$}}

\newcommand\ot{\otimes}

\newcommand\veps{\varepsilon}
\newcommand\vphi{\varphi}
\newcommand\wt{\mbox{\sl wt}\,}
\newcommand\Z{{\mathbb Z}}
\newcommand\Zn{\Z_{\ge0}}

\begin{document}

\allowdisplaybreaks

\renewcommand{\thefootnote}{$\star$}

\renewcommand{\PaperNumber}{022}

\FirstPageHeading

\ShortArticleName{Zero Action on Perfect Crystals for $U_q\big(G_2^{(1)}\big)$}

\ArticleName{Zero Action on Perfect Crystals for $\boldsymbol{U_q\big(G_2^{(1)}\big)}$\footnote{This paper is a
contribution to the Proceedings of the Workshop ``Geometric Aspects of Discrete and Ultra-Discrete Integrable Systems'' (March 30 -- April 3, 2009, University of Glasgow, UK). The full collection is
available at
\href{http://www.emis.de/journals/SIGMA/GADUDIS2009.html}{http://www.emis.de/journals/SIGMA/GADUDIS2009.html}}}

\Author{Kailash C. MISRA~$^\dag$, Mahathir MOHAMAD~$^\ddag$ and Masato OKADO~$^\ddag$}

\AuthorNameForHeading{K.C. Misra, M. Mohamad and M. Okado}

\Address{$^\dag$~Department of Mathematics, North Carolina State University,\\
\hphantom{$^\dag$}~Raleigh, North Carolina 27695-8205, USA}
\EmailD{\href{mailto:misra@unity.ncsu.edu}{misra@unity.ncsu.edu}}

\Address{$^\ddag$~Department of Mathematical Science, Graduate School of Engineering Science,\\
\hphantom{$^\ddag$}~Osaka University, Toyonaka, Osaka 560-8531, Japan}
\EmailD{\href{mailto:mahathir75@yahoo.com}{mahathir75@yahoo.com}, \href{mailto:okado@sigmath.es.osaka-u.ac.jp}{okado@sigmath.es.osaka-u.ac.jp}}

\ArticleDates{Received November 13, 2009, in f\/inal form March 03, 2010;  Published online March 09, 2010}

\Abstract{The actions of $0$-Kashiwara operators on the $U'_q\big(G_2^{(1)}\big)$-crystal $B_l$
in [Yamane S., {\it J.~Algebra} {\bf 210} (1998), 440--486] are made explicit by using a similarity technique from that of a $U'_q\big(D_4^{(3)}\big)$-crystal.
It is shown that $\{B_l\}_{l\ge1}$ forms a coherent family of perfect crystals.}

\Keywords{combinatorial representation theory; quantum af\/f\/ine algebra; crystal bases}

\Classification{05E99; 17B37; 17B67; 81R10; 81R50}

\section{Introduction}

Let $\geh$ be a symmetrizable Kac--Moody algebra. Let $I$ be its index set for simple roots,
$P$ the weight lattice, $\alpha_i\in P$ a simple root ($i\in I$), and $h_i\in P^*(=\mbox{Hom}
(P,\Z))$ a simple coroot ($i\in I$). To each $i\in I$ we associate a positive integer $m_i$
and set $\tilde{\alpha}_i=m_i\alpha_i$, $\tilde{h}_i=h_i/m_i$. Suppose $(\langle\tilde{h}_i,
\tilde{\alpha}_j\rangle)_{i,j\in I}$ is a generalized Cartan matrix for another symmetrizable
Kac--Moody algebra $\tilde{\geh}$. Then the subset $\tilde{P}$ of $P$ consisting of $\la\in P$
such that $\langle\tilde{h}_i,\la\rangle$ is an integer for any $i\in I$ can be considered as
the weight lattice of $\tilde{\geh}$. For a dominant integral weight $\la$ let $B^{\geh}(\la)$
be the highest weight crystal with highest weight $\la$ over $U_q(\geh)$. Then, in \cite{K2.5}
Kashiwara showed the following. (The theorem in \cite{K2.5} is more general.)

\begin{theorem} \label{th:Kashiwara}
Let $\la$ be a dominant integral weight in $\tilde{P}$. Then, there exists a unique injective map
$S:B^{\tilde{\geh}}(\la)\rightarrow B^{\geh}(\la)$ such that
\begin{gather*}
\wt S(b) =\wt b,\qquad
S(e_ib) =e_i^{m_i}S(b),\qquad
S(f_ib) =f_i^{m_i}S(b).
\end{gather*}
\end{theorem}

In this paper, we use this theorem to examine the so-called Kirillov--Reshetikhin crystal. Let~$\geh$ be the af\/f\/ine algebra of type $D_4^{(3)}$. The generalized Cartan matrix
$(\langle h_i,\alpha_j\rangle)_{i,j\in I}$ ($I=\{0,1,2\}$) is given by
\[
\left(
\begin{array}{rrr}
2  & -1 & 0  \\
-1 & 2  & -3 \\
0  & -1 & 2
\end{array}
\right).
\]
Set $(m_0,m_1,m_2)=(3,3,1)$. Then, $\tilde{\geh}$ def\/ined above turns out to be the af\/f\/ine algebra of type~$G_2^{(1)}$. Their Dynkin diagrams are depicted as follows
\[
\lower1mm\hbox{$D_4^{(3)}:$}\quad
\begin{minipage}{35mm}
\xymatrix@R=.0ex{
0&1&2\\
\circ\ar@{-}[r]&\circ&\circ\ar@3{->}[l]
}
\end{minipage}\qquad\qquad
\lower1mm\hbox{$G_2^{(1)}:$}\quad
\begin{minipage}{35mm}
\xymatrix@R=.0ex{
0&1&2\\
\circ\ar@{-}[r]&\circ&\circ\ar@3{<-}[l]
}
\end{minipage}
\]
For $G_2^{(1)}$ a family of perfect crystals $\{B_l\}_{l\ge1}$ was constructed in \cite{Y}. However,
the crystal elements there were realized in terms of tableaux given in \cite{KM}, and it was not easy
to calculate the action of $0$-Kashiwara operators on these tableaux. On the other hand,
an explicit action of these operators was given on perfect crystals $\{\hat{B}_l\}_{l\ge1}$ over
$U'_q\big(D_4^{(3)}\big)$ in \cite{KMOY}. Hence, it is a natural idea to use Theorem \ref{th:Kashiwara} to
obtain the explicit action of $e_0,f_0$ on $B_l$ from that on $\hat{B}_{l'}$ with suitable~$l'$.
We remark that Kirillov--Reshetikhin crystals are parametrized by a node of the Dynkin diagram except
$0$ and a positive integer. Both $B_l$ and $\hat{B}_l$ correspond to the pair~$(1,l)$.

Our strategy to do this is as follows. We def\/ine $V_l$ as an appropriate subset of $\hat{B}_{3l}$ that
is closed under the action of $\hat{e}_i^{m_i}$, $\hat{f}_i^{m_i}$ where $\hat{e}_i$, $\hat{f}_i$ stand for
the Kashiwara operators on~$\hat{B}_{3l}$. Hence, we can regard $V_l$ as a $U'_q\big(G_2^{(1)}\big)$-crystal.
We next show that as a $U_q\big(G_2^{(1)}\big)_{\{0,1\}}(=U_q(A_2))$-crystal
and as a $U_q\big(G_2^{(1)}\big)_{\{1,2\}}(=U_q(G_2))$-crystal, $V_l$ has the same decomposition as~$B_l$.
Then, we can conclude from Theorem~6.1 of~\cite{KMOY} that~$V_l$ is isomorphic to the
$U'_q\big(G_2^{(1)}\big)$-crystal $B_l$ constructed in~\cite{Y} (Theorem~\ref{th:main}).

The paper is organized as follows. In Section~\ref{section2} we review the $U'_q\big(D_4^{(3)}\big)$-crystal $\hat{B}_l$.
We then construct a $U'_q\big(G_2^{(1)}\big)$-crystal $V_l$ in $\hat{B}_{3l}$ with the aid of Theorem~\ref{th:Kashiwara} and see it coincides with $B_l$ given in \cite{Y} in Section~\ref{sec:G2(1)-crystal}. Minimal elements
of $B_l$ are found and $\{B_l\}_{l\ge1}$ is shown to form a~coherent family of perfect crystals in
Section~\ref{section4}. The crystal graphs of $B_1$ and $B_2$ are included in Section~\ref{section5}.

\section[Review on $U'_q(D_4^{(3)})$-crystal $\hat{B}_l$]{Review on $\boldsymbol{U'_q\big(D_4^{(3)}\big)}$-crystal $\boldsymbol{\hat{B}_l}$}\label{section2}

In this section we recall the perfect crystal for $U'_q\big(D_4^{(3)}\big)$ constructed in \cite{KMOY}.
Since we also consider $U'_q\big(G_2^{(1)}\big)$-crystals later, we denote it by $\hat{B}_l$. Kashiwara
operators~$e_i$, $f_i$ and $\veps_i$, $\vphi_i$ on~$\hat{B}_l$ are denoted by $\hat{e}_i$, $\hat{f}_i$ and
$\hat{\veps}_i$, $\hat{\vphi}_i$. Readers are warned that the coordinates $x_i$, ${\bar x}_i$
and steps by Kashiwara operators in~\cite{KMOY} are divided by~3 here, since it is more convenient
for our purpose. As a set
\[
\hat{B}_l=\left\{b=(x_1,x_2,x_3,{\bar x}_3,{\bar x}_2,{\bar x}_1)\in(\Zn/3)^6
\left\vert
\begin{array}{l}
3x_3\equiv3\bar{x}_3\;(\text{mod }2), \vspace{1mm}\\
{\displaystyle \sum_{i=1,2}} (x_i+{\bar x}_i)+(x_3+{\bar x}_3)/2 \leq l/3
\end{array}
\right.
\right\}.
\]

In order to def\/ine the actions of Kashiwara operators $\hat{e}_i$ and
$\hat{f}_i$ for $i= 0,1,2$, we introduce some notations and conditions.
Set $(x)_+=\max(x,0)$. For $b=(x_1,x_2,x_3,\bar{x}_3,\bar{x}_2,\bar{x}_1)\in \hat{B}_l$ we set
\begin{equation} \label{def s(b)}
s(b)=x_1+x_2+\frac{x_3+{\bar x}_3}{2}+{\bar x}_2+{\bar x}_1,
\end{equation}
and
\begin{equation} \label{z1-4}
z_1={\bar x}_1-x_1, \qquad
z_2={\bar x}_2 -{\bar x}_3, \qquad
z_3=x_3-x_2, \qquad
z_4=({\bar x}_3-x_3)/2.
\end{equation}

Now we def\/ine conditions ($E_1$)--($E_6$) and ($F_1$)--($F_6$) as follows
\begin{equation} \label{(F)}
\begin{split}
& (F_1)\quad
z_1+z_2+z_3+3z_4\le0,\quad z_1+z_2+3z_4\le0, \quad z_1+z_2\le0, \quad z_1\le0,
\\
& (F_2)\quad
z_1+z_2+z_3+3z_4\le0, \quad z_2+3z_4\le0, \quad z_2\le0, \quad z_1> 0,
\\
& (F_3)\quad
z_1+z_3+3z_4\le0, \quad z_3+3z_4\le0, \quad z_4\le0, \quad z_2> 0, \quad z_1+z_2> 0,
\\
& (F_4)\quad
z_1+z_2+3z_4> 0, \quad z_2+3z_4> 0, \quad z_4> 0, \quad z_3\le0, \quad z_1+z_3\le0,
\\
& (F_5)\quad
z_1+z_2+z_3+3z_4> 0, \quad z_3+3z_4> 0, \quad z_3> 0, \quad z_1\le0,
\\
& (F_6)\quad
z_1+z_2+z_3+3z_4> 0, \quad z_1+z_3+3z_4> 0,\quad z_1+z_3> 0,\quad z_1> 0.
\end{split}
\end{equation}
The conditions ($F_1$)--($F_6$) are disjoint and they exhaust all cases.
($E_i$) ($1\le i\le 6$) is def\/ined from ($F_i$) by replacing $>$ (resp. $\le$) with $\ge$ (resp. $<$).
We also def\/ine
\begin{equation} \label{A}
A=(0,z_1,z_1+z_2,z_1+z_2+3z_4,z_1+z_2+z_3+3z_4,2z_1+z_2+z_3+3z_4).
\end{equation}
Then, for $b=(x_1,x_2,x_3,{\bar x}_3,{\bar x}_2,{\bar x}_1) \in \hat{B}_l$, $\hat{e}_ib$,
$\hat{f}_ib$, $\hat{\veps}_i(b)$, $\hat{\vphi}_i(b)$ are given as follows
\begin{gather}
\hat{e}_0 b=
\begin{cases}
(x_1 -1/3,\ldots)
& \text{if ($E_1$)},
\\
(\ldots,x_3 -1/3,{\bar x}_3 -1/3,\ldots,{\bar x}_1 +1/3)
& \text{if ($E_2$)},
\\
(\ldots,x_3 -2/3,\ldots,{\bar x}_2 +1/3,\ldots)
& \text{if ($E_3$)},
\\
(\ldots,x_2 -1/3,\ldots,{\bar x}_3 +2/3,\ldots)
& \text{if ($E_4$)},
\\
(x_1 -1/3,\ldots,x_3 +1/3,{\bar x}_3 +1/3,\ldots)
& \text{if ($E_5$)},
\\
(\ldots,{\bar x}_1 +1/3)
& \text{if ($E_6$)},
\end{cases}
\nonumber\\
\hat{f}_0 b=
\begin{cases}
(x_1 +1/3,\ldots)
& \text{if ($F_1$)},
\\
(\ldots,x_3 +1/3,{\bar x}_3 +1/3,\ldots,{\bar x}_1 -1/3)
& \text{if ($F_2$)},
\\
(\ldots,x_3 +2/3,\ldots,{\bar x}_2 -1/3,\ldots)
& \text{if ($F_3$)},
\\
(\ldots,x_2 +1/3,\ldots,{\bar x}_3 -2/3,\ldots)
& \text{if ($F_4$)},
\\
(x_1 +1/3,\ldots,x_3 -1/3,{\bar x}_3 -1/3,\ldots)
& \text{if ($F_5$)},
\\
(\ldots,{\bar x}_1 -1/3)
& \text{if ($F_6$)}.
\end{cases}
\nonumber\\
\hat{e}_1 b=
\begin{cases}
(\ldots,{\bar x}_2 +1/3,{\bar x}_1 -1/3)
& \text{if $z_2 \geq (-z_3)_+$},
\\
(\ldots,x_3 +1/3,{\bar x}_3 -1/3,\ldots)
& \text{if $z_2 <0\leq z_3$},
\\
(x_1 +1/3,x_2 -1/3,\ldots)
& \text{if $(z_2)_+ < (-z_3)$},
\end{cases}
\nonumber\\
\hat{f}_1 b=
\begin{cases}
(x_1 -1/3,x_2 +1/3,\ldots)
& \text{if $(z_2)_+ \leq (-z_3)$},
\\
(\ldots,x_3 -1/3,{\bar x}_3 +1/3,\ldots)
& \text{if $z_2 \leq 0< z_3$},
\\
(\ldots,{\bar x}_2 -1/3,{\bar x}_1 +1/3)
& \text {if $ z_2 >(-z_3)_+$},
\end{cases}
\nonumber\\
\hat{e}_2 b=
\begin{cases}
(\ldots,{\bar x}_3 +2/3,{\bar x}_2 -1/3,\ldots)
& \text{if $z_4 \geq 0$},
\\
(\ldots,x_2 +1/3,x_3 -2/3,\ldots)
& \text{if $z_4 < 0$},
\end{cases}
\nonumber\\
\hat{f}_2 b=
\begin{cases}
(\ldots,x_2 -1/3,x_3 +2/3,\ldots)
& \text{if $z_4 \leq 0$},
\\
(\ldots,{\bar x}_3 -2/3,{\bar x}_2 +1/3,\ldots)
& \text{if $z_4 > 0$},
\end{cases}\nonumber
\\
 \label{eps-phi hat}
\begin{split}
& \hat{\veps}_0(b) =l-3s(b)+3\max A-3(2z_1+z_2+z_3+3z_4),\\
& \hat{\vphi}_0(b) =l-3s(b)+3\max A,\\
& \hat{\veps}_1(b) =3{\bar x}_1+3({\bar x}_3-{\bar x}_2+(x_2-x_3)_+)_+,\\
& \hat{\vphi}_1(b) =3x_1+3(x_3-x_2+({\bar x}_2-{\bar x}_3)_+)_+,\\
& \hat{\veps}_2(b) =3{\bar x}_2+\tfrac{3}{2}(x_3-{\bar x}_3)_+,\qquad
  \hat{\vphi}_2(b)=3x_2+\tfrac{3}{2}({\bar x}_3-x_3)_+.
\end{split}
\end{gather}
If $\hat{e}_ib$ or $\hat{f}_ib$ does not belong to $\hat{B}_l$, namely, if $x_j$ or $\bar{x}_j$
for some $j$ becomes negative or $s(b)$ exceeds $l/3$, we should understand it to be $0$.
Forgetting the $0$-arrows,
\[
\hat{B}_l\simeq\bigoplus_{j=0}^lB^{G^\dagger_2}(j\La_1),
\]
where $B^{G^\dagger_2}(\la)$ is the highest weight $U_q(G^\dagger_2)$-crystal of highest weight
$\la$ and $G^\dagger_2$ stands for the simple Lie algebra $G_2$ with the reverse labeling of
the indices of the simple roots ($\alpha_1$ is the short root).
Forgetting $2$-arrows,
\[
\hat{B}_l\simeq
\bigoplus^{\lfloor\frac{l}{2}\rfloor}_{i=0}
\bigoplus_{\begin{subarray}{c}
{i\leq j_0, j_1 \leq l-i}\\
{j_0, j_1 \equiv l-i  \!\!\!\!\pmod 3}
\end{subarray}}
B^{A_2}(j_0\La_0+j_1\La_1),
\]
where $B^{A_2}(\la)$ is the highest weight $U_q(A_2)$-crystal (with indices $\{0,1\}$)
of highest weight $\la$.

\section[$U'_q(G_2^{(1)})$-crystal]{$\boldsymbol{U'_q\big(G_2^{(1)}\big)}$-crystal} \label{sec:G2(1)-crystal}

In this section we def\/ine a subset $V_l$ of $\hat{B}_{3l}$ and see it is isomorphic to the
$U'_q\big(G_2^{(1)}\big)$-crystal $B_l$. The set $V_l$ is def\/ined as a subset of $\hat{B}_{3l}$
satisfying the following conditions:
\begin{equation} \label{constr}
x_1,{\bar x}_1,x_2-x_3,{\bar x}_3-{\bar x}_2\in\Z.
\end{equation}
For an element $b=(x_1,x_2,x_3,{\bar x}_3,{\bar x}_2,{\bar x}_1)$ of $V_l$ we def\/ine $s(b)$ as in~\eqref{def s(b)}.
From \eqref{constr} we see that $s(b)\in\{0,1,\ldots,l\}$.

\begin{lemma} \label{lem:count}
For $0\le k\le l$
\[
\sharp\{b\in V_l\mid s(b)=k\}=\frac1{120}(k+1)(k+2)(2k+3)(3k+4)(3k+5).
\]
\end{lemma}

\begin{proof}
We f\/irst count the number of elements $(x_2,x_3,{\bar x}_3,{\bar x}_2)$ satisfying the conditions
of coordinates as an element of $V_l$ and $x_2+(x_3+{\bar x}_3)/2+{\bar x}_2=m$ ($m=0,1,\ldots,k$).
According to $(a,b,c,d)$ ($a,d\in\{0,1/3,2/3\}$, $b,c\in\{0,1/3,2/3,1,4/3,5/3\}$) such that
$x_2\in\Z+a$, $x_3\in2\Z+b$, ${\bar x}_3\in2\Z+c$, ${\bar x}_2\in\Z+d$,
we divide the cases into the following 18:
\begin{alignat*}{3}
&\mbox{(i)}\ (0,0,0,0),&&\mbox{(ii)}\ (0,0,2/3,2/3),&&\mbox{(iii)}\ (0,0,4/3,1/3),\\
&\mbox{(iv)}\ (0,1,1/3,1/3),&&\mbox{(v)}\ (0,1,1,0),&&\mbox{(vi)}\ (0,1,5/3,2/3),\\
&\mbox{(vii)}\ (1/3,1/3,1/3,1/3),&&\mbox{(viii)}\ (1/3,1/3,1,0),&&\mbox{(ix)}\ (1/3,1/3,5/3,2/3),\\
&\mbox{(x)}\ (1/3,4/3,0,0),&&\mbox{(xi)}\ (1/3,4/3,2/3,2/3),&&\mbox{(xii)}\ (1/3,4/3,4/3,1/3),\\
&\mbox{(xiii)}\ (2/3,2/3,0,0),&& \mbox{(xiv)}\ (2/3,2/3,2/3,2/3),\qquad &&\mbox{(xv)}\ (2/3,2/3,4/3,1/3),\\
&\mbox{(xvi)}\ (2/3,5/3,1/3,1/3), \qquad && \mbox{(xvii)}\ (2/3,5/3,1,0),\qquad & & \mbox{(xviii)}\ (2/3,5/3,5/3,2/3).
\end{alignat*}
The number of elements $(x_2,x_3,{\bar x}_3,{\bar x}_2)$ in a case among the above such that
$a+(b+c)/2+d=e$ ($e=0,1,2,3$) is given by $f(e)={m-e+3\choose3}$.
Since there is one case with $e=0$ (i) and $e=3$ (xviii) and 8 cases with $e=1$ and $e=2$,
the number of $(x_2,x_3,{\bar x}_3,{\bar x}_2)$ such that $x_2+(x_3+{\bar x}_3)/2+{\bar x}_2=m$
is given by
\[
f(0)+8f(1)+8f(2)+f(3)=\frac12(2m+1)(3m^2+3m+2).
\]
For each $(x_2,x_3,{\bar x}_3,{\bar x}_2)$ such that $x_2+(x_3+{\bar x}_3)/2+{\bar x}_2=m$
($m=0,1,\ldots,k$) there are $(k-m+1)$ cases for $(x_1,{\bar x}_1)$, so the number of
$b\in V_l$ such that $s(b)=k$ is given by
\[
\sum_{m=0}^k\frac12(2m+1)(3m^2+3m+2)(k-m+1).
\]
A direct calculation leads to the desired result.
\end{proof}

We def\/ine the action of operators $e_i,f_i$ ($i=0,1,2$) on $V_l$ as follows.
\begin{gather*}
e_0 b=
\begin{cases}
(x_1 -1,\ldots)
& \text{if ($E_1$)},
\\
(\ldots,x_3 -1,{\bar x}_3 -1,\ldots,{\bar x}_1 +1)
& \text{if ($E_2$)},
\\
\big(\ldots,x_2-\frac23,x_3-\frac23,{\bar x}_3+\frac43,{\bar x}_2+\frac13,\ldots\big)
& \text{if ($E_3$) and $z_4=-\frac13$},
\\
\big(\ldots,x_2-\frac13,x_3-\frac43,{\bar x}_3+\frac23,{\bar x}_2+\frac23,\ldots\big)
& \text{if ($E_3$) and $z_4=-\frac23$},
\\
(\ldots,x_3 -2,\ldots,{\bar x}_2 +1,\ldots)
& \text{if ($E_3$) and $z_4\ne-\frac13,-\frac23$},
\\
(\ldots,x_2 -1,\ldots,{\bar x}_3 +2,\ldots)
& \text{if ($E_4$)},
\\
(x_1 -1,\ldots,x_3 +1,{\bar x}_3 +1,\ldots)
& \text{if ($E_5$)},
\\
(\ldots,{\bar x}_1 +1)
& \text{if ($E_6$)},
\end{cases}
\\
f_0 b=
\begin{cases}
(x_1 +1,\ldots)
& \text{if ($F_1$)},
\\
(\ldots,x_3 +1,{\bar x}_3 +1,\ldots,{\bar x}_1 -1)
& \text{if ($F_2$)},
\\
(\ldots,x_3 +2,\ldots,{\bar x}_2 -1,\ldots)
& \text{if ($F_3$)},
\\
\big(\ldots,x_2 +\frac13,x_3+\frac43,{\bar x}_3-\frac23,{\bar x}_2 -\frac23,\ldots\big)
& \text{if ($F_4$) and $z_4=\frac13$},
\\
\big(\ldots,x_2 +\frac23,x_3+\frac23,{\bar x}_3-\frac43,{\bar x}_2 -\frac13,\ldots\big)
& \text{if ($F_4$) and $z_4=\frac23$},
\\
(\ldots,x_2 +1,\ldots,{\bar x}_3 -2,\ldots)
& \text{if ($F_4$) and $z_4\ne\frac13,\frac23$},
\\
(x_1 +1,\ldots,x_3 -1,{\bar x}_3 -1,\ldots)
& \text{if ($F_5$)},
\\
(\ldots,{\bar x}_1 -1)
& \text{if ($F_6$)},
\end{cases}
\\
e_1 b=
\begin{cases}
(\ldots,{\bar x}_2 +1,{\bar x}_1 -1)
& \text{if ${\bar x}_2 -{\bar x}_3 \geq (x_2 -x_3)_+$},
\\
(\ldots,x_3 +1,{\bar x}_3 -1,\ldots)
& \text{if ${\bar x}_2 -{\bar x}_3 <0\leq x_3 -x_2$},
\\
(x_1 +1,x_2 -1,\ldots)
& \text{if $({\bar x}_2 -{\bar x}_3)_+ <x_2 -x_3$},
\end{cases}
\\
f_1 b=
\begin{cases}
(x_1 -1,x_2 +1,\ldots)
& \text{if $({\bar x}_2 -{\bar x}_3)_+ \leq x_2 -x_3$},
\\
(\ldots,x_3 -1,{\bar x}_3 +1,\ldots)
& \text{if ${\bar x}_2 -{\bar x}_3 \leq 0<x_3 -x_2$},
\\
(\ldots,{\bar x}_2 -1,{\bar x}_1 +1)
& \text {if ${\bar x}_2 -{\bar x}_3 >(x_2 -x_3)_+$},
\end{cases}
\\
e_2 b=
\begin{cases}
\big(\ldots,{\bar x}_3 +\frac23,{\bar x}_2 -\frac13,\ldots\big)
& \text{if ${\bar x}_3 \geq x_3$},
\\
\big(\ldots,x_2 +\frac13,x_3 -\frac23,\ldots\big)
& \text{if ${\bar x}_3 <x_3$},
\end{cases}
\\
f_2 b=
\begin{cases}
\big(\ldots,x_2 -\frac13,x_3 +\frac23,\ldots\big)
& \text{if ${\bar x}_3 \leq x_3$},
\\
\big(\ldots,{\bar x}_3 -\frac23,{\bar x}_2 +\frac13,\ldots\big)
& \text{if ${\bar x}_3 >x_3$}.
\end{cases}
\end{gather*}

We now set $(m_0,m_1,m_2)=(3,3,1)$.

\begin{proposition} \label{prop:imp}\qquad

\begin{itemize}\itemsep=0pt
\item[$(1)$] For any $b\in V_l$ we have $e_ib,f_ib\in V_l\sqcup\{0\}$ for $i=0,1,2$.
\item[$(2)$] The equalities $e_i=\hat{e}_i^{m_i}$ and $f_i=\hat{f}_i^{m_i}$
	hold on $V_l$ for $i=0,1,2$.
\end{itemize}
\end{proposition}

\begin{proof}
(1) can be checked easily.

For (2) we only treat $f_i$. To prove the $i=0$ case consider the following table
\begin{center}
\begin{tabular}{c|cccccc}
&($F_1$)&($F_2$)&($F_3$)&($F_4$)&($F_5$)&($F_6$)\\
\hline
$z_1$&$-1/3$&$-1/3$&0&0&$-1/3$&$-1/3$\\
$z_2$&0&$-1/3$&$-1/3$&2/3&1/3&0\\
$z_3$&0&1/3&2/3&$-1/3$&$-1/3$&0\\
$z_4$&0&0&$-1/3$&$-1/3$&0&0
\end{tabular}
\end{center}
This table signif\/ies the dif\/ference $(z_j\text{ for }\hat{f}_0b)-(z_j\text{ for }b)$ when
$b$ belongs to the case $(F_i)$. Note that the left hand sides of the inequalities of each
$(F_i)$ \eqref{(F)} always decrease by $1/3$. Since $z_1,z_2,z_3\in\Z,z_4\in\Z/3$ for $b\in V_l$,
we see that if $b$ belongs to $(F_i)$, $\hat{f}_0b$ and $\hat{f}_0^2b$ also belong to $(F_i)$
except two cases: (a) $b\in(F_4)$ and $z_4=1/3$, and (b) $b\in(F_4)$ and $z_4=2/3$. If (a)
occurs, we have $\hat{f}_0b,\hat{f}_0^2b\in(F_3)$. Hence, we obtain $f_0=\hat{f}_0^3$ in
this case. If (b) occurs, we have $\hat{f}_0b\in(F_4)$, $\hat{f}_0^2b\in(F_3)$.
Therefore, we obtain $f_0=\hat{f}_0^3$ in this case as well.

In the $i=1$ case, if $b$ belongs to one of the 3 cases, $\hat{f}_1b$ and $\hat{f}_1^2b$ also
belong to the same case. Hence, we obtain $f_1=\hat{f}_1^3$. For $i=2$ there is nothing to do.
\end{proof}

Proposition \ref{prop:imp}, together with Theorem \ref{th:Kashiwara}, shows that $V_l$ can be
regarded as a $U'_q\big(G_2^{(1)}\big)$-crystal with operators $e_i$, $f_i$ ($i=0,1,2$).

\begin{proposition} \label{prop:G2}
As a $U_q\big(G_2^{(1)}\big)_{\{1,2\}}$-crystal
\begin{gather*}
V_l\simeq \bigoplus_{k=0}^lB^{G_2}(k\La_1),
\end{gather*}
where $B^{G_2}(\la)$ is the highest weight $U_q(G_2)$-crystal of highest weight~$\la$.
\end{proposition}

\begin{proof}
For a subset $J$ of $\{0,1,2\}$ we say $b$ is $J$-highest if $e_jb=0$ for any $j\in J$.
Note from~\eqref{eps-phi hat} that $b_k=(k,0,0,0,0,0)$ ($0\le k\le l$) is $\{1,2\}$-highest of
weight $3k\La_1$ in $\hat{B}_{3l}$. By setting $\geh=G^\dagger_2$ ($=G_2$ with the reverse labeling
of indices), $(m_1,m_2)=(3,1)$, $\tilde{\geh}=G_2$ in Theorem~\ref{th:Kashiwara}, we know that the
connected component generated from $b_k$ by $f_1=\hat{f}_1^3$ and $f_2=\hat{f}_2$ is
isomorphic to $B^{G_2}(k\La_1)$. Hence by Proposition~\ref{prop:imp}~(1) we have
\begin{equation} \label{incl}
\bigoplus_{k=0}^lB^{G_2}(k\La_1)\subset V_l.
\end{equation}
Now recall Weyl's formula to calculate the dimension of the highest weight representation.
In our case we obtain
\[
\sharp B^{G_2}(k\La_1)=\frac1{120}(k+1)(k+2)(2k+3)(3k+4)(3k+5).
\]
However, this is equal to $\sharp\{b\in V_l\mid s(b)=k\}$ by Lemma \ref{lem:count}. Therefore,
$\subset$ in \eqref{incl} should be $=$, and the proof is completed.
\end{proof}

\begin{proposition}
As a $U_q\big(G_2^{(1)}\big)_{\{0,1\}}$-crystal
\[
V_l\simeq \bigoplus_{i=0}^{\lfloor l/2\rfloor}\bigoplus_{i\le j_0,j_1\le l-i}
B^{A_2}(j_0\La_0+j_1\La_1),
\]
where $B^{A_2}(\la)$ is the highest weight $U_q(A_2)$-crystal $($with indices $\{0,1\})$
of highest weight~$\la$.
\end{proposition}

\begin{proof}
For integers $i,j_0,j_1$ such that $0\le i\le l/2$, $i\le j_0,j_1\le l-i$, def\/ine an element
$b_{i,j_0,j_1}$ of~$V_l$ by
\[
b_{i,j_0,j_1}=\left\{
\begin{array}{ll}
(0,y_1,3y_0-2y_1+i,y_0+i,y_0+j_0,0)\quad&\text{if }j_0\le j_1,\\
(0,y_0,y_0+i,2y_1-y_0+i,2y_0-y_1+j_0,0)\quad&\text{if }j_0>j_1.
\end{array}\right.
\]
Here we have set $y_a=(l-i-j_a)/3$ for $a=0,1$. From \eqref{eps-phi hat} one notices that $b_{i,j_0,j_1}$
is $\{0,1\}$-highest of weight $3j_0\La_0+3j_1\La_1$ in $\hat{B}_{3l}$. For instance,
$\hat{\veps}_0(b_{i,j_0,j_1})=0$ and $\hat{\vphi}_0(b_{i,j_0,j_1})=3j_0$ since $s(b_{i,j_0,j_1})=l$
and $\max A=2z_1+z_2+z_3+3z_4=j_0$.
By setting $\geh=\tilde{\geh}=A_2,
(m_0,m_1)=(3,3)$ in Theorem~\ref{th:Kashiwara}, the connected component generated from $b_{i,j_0,j_1}$ by
$f_i=\hat{f}_i^3$ ($i=0,1$) is isomorphic to $B^{A_2}(j_0\La_0+j_1\La_1)$. Hence, by Proposition~\ref{prop:imp}~(1) we have
\[
\bigoplus_{i=0}^{\lfloor l/2\rfloor}\bigoplus_{i\le j_0,j_1\le l-i}
B^{A_2}(j_0\La_0+j_1\La_1)\subset V_l.
\]
However, from Proposition \ref{prop:G2} one knows that
\[
\sharp V_l=\sum_{k=0}^l\sharp B^{G_2}(k\La_1).
\]
Moreover, it is already established in \cite{Y} that
\[
\sum_{k=0}^l\sharp B^{G_2}(k\La_1)=\sum_{i=0}^{\lfloor l/2\rfloor}\sum_{i\le j_0,j_1\le l-i}
\sharp B^{A_2}(j_0\La_0+j_1\La_1).
\]
Therefore, the proof is completed.
\end{proof}

Theorem 6.1 in \cite{KMOY} shows that if two $U'_q\big(G_2^{(1)}\big)$-crystals decompose into
$\bigoplus_{0\le k\le l}B^{G_2}(k\La_1)$ as $U_q(G_2)$-crystals, then they are isomorphic to
each other. Therefore, we now have

\begin{theorem} \label{th:main}
$V_l$ agrees with the $U'_q\big(G_2^{(1)}\big)$-crystal $B_l$ constructed in {\rm \cite{Y}}.
The values of $\veps_i$, $\vphi_i$ with our representation are given by
\begin{equation} \label{eps-phi}
\begin{split}
& \veps_0(b) =l-s(b)+\max A-(2z_1+z_2+z_3+3z_4),\qquad \vphi_0(b) =l-s(b)+\max A,\\
& \veps_1(b) ={\bar x}_1+({\bar x}_3-{\bar x}_2+(x_2-x_3)_+)_+,\qquad
 \varphi_1(b) =x_1+(x_3-x_2+({\bar x}_2-{\bar x}_3)_+)_+,\\
& \veps_2(b) =3{\bar x}_2+\tfrac{3}{2}(x_3-{\bar x}_3)_+,\qquad
\varphi_2(b)=3x_2+\tfrac{3}{2}({\bar x}_3-x_3)_+.
\end{split}
\end{equation}
\end{theorem}

\section{Minimal elements and a coherent family}\label{section4}

The notion of perfect crystals was introduced in \cite{KMN1} to construct the path realization
of a highest weight crystal of a quantum af\/f\/ine algebra. The crystal $B_l$ was shown to be
perfect of level~$l$ in~\cite{Y}. In this section we obtain all the minimal
elements of~$B_l$ in the coordinate representation and also show $\{B_l\}_{l\ge1}$ forms a coherent
family of perfect crystals. For the notations such as~$P_{cl}$, $(P_{cl}^+)_l$, see~\cite{KMN1}.

\subsection{Minimal elements}

From \eqref{eps-phi} we have
\begin{gather*}
\langle c, \varphi(b)\rangle
= \varphi_0(b)+2\varphi_1(b)+\varphi_2(b) \\
\phantom{\langle c, \varphi(b)\rangle}{} = l+\max A+2(z_3+(z_2)_+)_++(3z_4)_+ -(z_1+z_2+2z_3+3z_4),
\end{gather*}
where $z_j$ ($1\le j\le 4$) are given in~\eqref{z1-4} and $A$ is given in~\eqref{A}.
The following lemma was proven in~\cite{KMOY}, although $\mathbb{Z}$ is replaced with
$\mathbb{Z}/3$ here.

\begin{lemma}\label{lem:PerfectCrystal}
For $(z_1,z_2,z_3,z_4)\in (\mathbb{Z}/3)^4$ set
\[
\psi(z_1,z_2,z_3,z_4)=\max A+
2(z_3+(z_2)_+)_+ +(3z_4)_+ -(z_1+z_2+2z_3+3z_4).
\]
Then we have $\psi(z_1,z_2,z_3,z_4)\geq 0$ and $\psi(z_1,z_2,z_3,z_4)=0$ if and only if
$(z_1,z_2,z_3,z_4)=(0,0,0,0)$.
\end{lemma}

From this lemma, we have $\langle c,\varphi(b)\rangle-l = \psi(z_1,z_2,z_3,z_4)\geq 0$.
Since $\langle c,\vphi(b)-\veps(b)\rangle=0$, we also have $\langle c,\veps(b)\rangle\ge l$.

Suppose $\langle c,\veps(b)\rangle=l$. It implies $\psi=0$. Hence from the lemma one can conclude
that such element $b=(x_1,x_2,x_3,\bar{x}_3,\bar{x}_2,\bar{x}_1)$ should satisfy
$x_1={\bar x}_1$, $x_2=x_3={\bar x}_3={\bar x}_2$. Therefore the set of minimal elements $(B_l)_{\min}$
in $B_l$ is given by
\[
(B_l)_{\min}=\{(\alpha,\beta,\beta,\beta,\beta,\alpha)
\vert\  \alpha\in \Zn,\beta\in(\Zn)/3,2\alpha+3\beta\le l\}.
\]
For $b=(\alpha,\beta,\beta,\beta,\beta,\alpha)\in (B_l)_{\min}$ one calculates
\[
\veps(b)=\vphi(b)=(l-2\alpha-3\beta)\La_0+\alpha\La_1+3\beta\La_2.
\]

\subsection{Coherent family of perfect crystals}

The notion of a coherent family of perfect crystals was introduced in~\cite{KKM}.
Let $\{B_l\}_{l\ge1}$ be a~family of perfect crystals $B_l$ of level $l$ and $(B_l)_{\min}$ be
the subset of minimal elements of $B_l$. Set $J=\{(l,b)\mid l\in\Z_{>0},b\in(B_l)_{\min}\}$.
Let $\sigma$ denote the isomorphism of $(P^+_{cl})_l$ def\/ined by $\sigma=\veps\circ\vphi^{-1}$.
For $\la\in P_{cl}$, $T_{\la}$ denotes a crystal with a unique element $t_{\la}$ def\/ined in \cite{K2}.
For our purpose the following facts are suf\/f\/icient. For any $P_{cl}$-weighted crystal
$B$ and $\la,\mu\in P_{cl}$ consider the crystal
\[
T_{\la}\ot B\ot T_{\mu}=\{t_{\la}\ot b\ot t_{\mu}\mid b\in B\}.
\]
The def\/inition of $T_\la$ and the tensor product rule of crystals imply
\begin{alignat*}{3}
& \et{i}(t_{\la}\ot b\ot t_{\mu}) =t_{\la}\ot\et{i}b\ot t_{\mu}, \qquad & &
\ft{i}(t_{\la}\ot b\ot t_{\mu})=t_{\la}\ot\ft{i}b\ot t_{\mu}, & \\
& \veps_i(t_{\la}\ot b\ot t_{\mu})=\veps_i(b)-\langle h_i,\la\rangle,\qquad  &&
  \vphi_i(t_{\la}\ot b\ot t_{\mu})=\vphi_i(b)+\langle h_i,\mu\rangle, &\\
& \wt(t_{\la}\ot b\ot t_{\mu})=\la+\mu+\wt b.&&&
\end{alignat*}

\begin{definition} \label{def:coherent}
A crystal $B_\infty$ with an element $b_\infty$ is called a limit of $\{B_l\}_{l\ge1}$ if it
satisf\/ies the following conditions:
\begin{itemize}\itemsep=0pt
\item[$\bullet$] $\wt b_\infty=0,\veps(b_\infty)=\vphi(b_\infty)=0$,

\item[$\bullet$] for any $(l,b)\in J$, there exists an embedding of crystals
\[
f_{(l,b)}:\ \ T_{\veps(b)}\ot B_l\ot T_{-\vphi(b)}\longrightarrow B_\infty
\]
sending $t_{\veps(b)}\ot b\ot t_{-\vphi(b)}$ to $b_\infty$,

\item[$\bullet$] $B_\infty=\bigcup_{(l,b)\in J}\mathrm{Im}\,f_{(l,b)}$.
\end{itemize}
If a limit exists for the family $\{B_l\}$, we say that $\{B_l\}$ is a coherent family of perfect
crystals.
\end{definition}

Let us now consider the following set
\[
B_\infty=\left\{b=(\nu_1,\nu_2,\nu_3,\bar{\nu}_3,\bar{\nu}_2,\bar{\nu}_1)\in(\Z/3)^6
\left\vert
\begin{array}{l}
\nu_1,\bar{\nu}_1,\nu_2-\nu_3,\bar{\nu}_3-\bar{\nu}_2\in\Z, \\
3\nu_3\equiv3\bar{\nu}_3\;(\text{mod }2)
\end{array}
\right.
\right\},
\]
and set $b_\infty=(0,0,0,0,0,0)$. We introduce the crystal structure on $B_\infty$ as follows.
The actions of~$e_i$,~$f_i$ ($i=0,1,2$) are def\/ined by the same rule as in Section~\ref{sec:G2(1)-crystal}
with $x_i$ and $\bar{x}_i$ replaced with~$\nu_i$ and~$\bar{\nu}_i$. The
only dif\/ference lies in the fact that $e_ib$ or $f_ib$ never becomes $0$, since we allow a~coordinate to be negative and there is no restriction for the sum $s(b)=\sum\limits_{i=1}^2(\nu_i+\bar{\nu}_i)
+(\nu_3+\bar{\nu}_3)/2$. For $\veps_i$, $\vphi_i$ with $i=1,2$ we adopt the formulas in Section~\ref{sec:G2(1)-crystal}. For $\veps_0$, $\vphi_0$ we def\/ine
\begin{gather*}
\veps_0(b)=-s(b)+\max A-(2z_1+z_2+z_3+3z_4),\qquad
\vphi_0(b) =-s(b)+\max A,
\end{gather*}
where $A$ is given in \eqref{A}
and $z_1$, $z_2$, $z_3$, $z_4$ are given in \eqref{z1-4} with $x_i$, $\bar{x}_i$ replaced by $\nu_i$, $\bar{\nu}_i$.
Note that $\wt b_\infty=0$ and $\veps_i(b_\infty)=\vphi_i(b_\infty)=0$ for $i=0,1,2$.

Let $b_0=(\alpha,\beta,\beta,\beta,\beta,\alpha)$ be an element of $(B_l)_{\min}$. Since
$\veps(b_0)=\vphi(b_0)$, one can set $\sigma=\mbox{id}$. Let $\la=\veps(b_0)$. For
$b=(x_1,x_2,x_3,\bar{x}_3,\bar{x}_2,\bar{x}_1)\in B_l$ we def\/ine a map
\[
f_{(l,b_0)}:\ \ T_{\la}\ot B_l\ot T_{-\la}\longrightarrow B_\infty
\]
by
\[
f_{(l,b_0)}(t_{\la}\ot b\ot t_{-\la})=b'=(\nu_1,\nu_2,\nu_3,\bar{\nu}_3,\bar{\nu}_2,\bar{\nu}_1),
\]
where
\begin{alignat*}{3}
& \nu_1=x_1-\alpha, \qquad & & \bar{\nu}_1=\bar{x}_1-\alpha, &\\
& \nu_j=x_j-\beta, \qquad & & \bar{\nu}_j=\bar{x}_j-\beta\quad (j=2,3).&
\end{alignat*}
Note that $s(b')=s(b)-(2\alpha+3\beta)$. Then we have
\begin{gather*}
\wt(t_{\la}\ot b\ot t_{-\la}) =\wt b=\wt b', \\
\vphi_0(t_{\la}\ot b\ot t_{-\la}) =\vphi_0(b)+\langle h_0,-\la\rangle \\
\phantom{\vphi_0(t_{\la}\ot b\ot t_{-\la})}{}
=\vphi_0(b')+(l-s(b))+s(b')-(l-2\alpha-3\beta)=\vphi_0(b'),\\
\vphi_1(t_{\la}\ot b\ot t_{-\la}) =\vphi_1(b)+\langle h_1,-\la\rangle=\vphi_1(b')+\alpha-\alpha
=\vphi_1(b'), \\
\vphi_2(t_{\la}\ot b\ot t_{-\la}) =\vphi_2(b)+\langle h_2,-\la\rangle=\vphi_2(b')+3\beta-3\beta
=\vphi_2(b').
\end{gather*}
$\veps_i(t_{\la}\ot b\ot t_{-\la})=\veps_i(b')$ ($i=0,1,2$) also follows from the above calculations.

{}From the fact that $(z_j\mbox{ for }b)=(z_j\mbox{ for }b')$ it is straightforward to check that
if $b,e_ib\in B_l$ (resp.~$b,f_ib\in B_l$), then
$f_{(l,b_0)}(e_i(t_{\la}\ot b\ot t_{-\la}))=e_if_{(l,b_0)}(t_{\la}\ot b\ot t_{-\la})$
(resp.~$f_{(l,b_0)}(f_i(t_{\la}\ot b\ot t_{-\la}))=f_if_{(l,b_0)}(t_{\la}\ot b\ot t_{-\la})$).
Hence $f_{(l,b_0)}$ is a crystal embedding. It is easy to see that
$f_{(l,b_0)}(t_{\la}\ot b_0\ot t_{-\la})=b_\infty$. We can also check
$B_\infty=\bigcup_{(l,b)\in J}\mathrm{Im}\,f_{(l,b)}$. Therefore we have shown that the family of
perfect crystals $\{B_l\}_{l\ge1}$ forms a coherent family.

\section[Crystal graphs of $B_1$ and $B_2$]{Crystal graphs of $\boldsymbol{B_1}$ and $\boldsymbol{B_2}$}
\label{section5}

\begin{figure}[t]

\centerline{\includegraphics{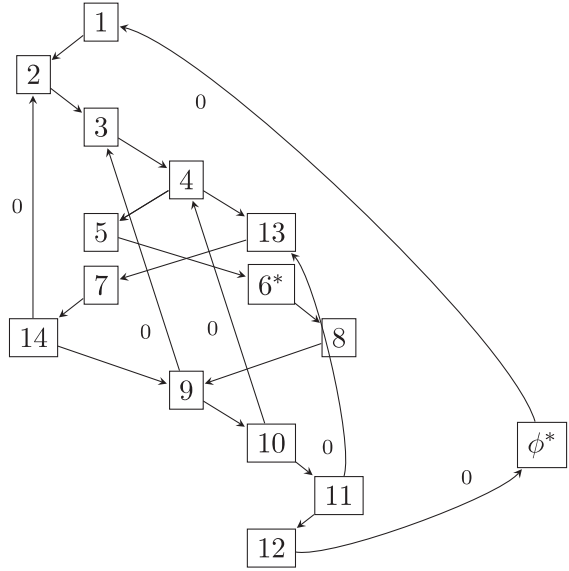}}

\caption{Crystal graph of $B_1$. $\swarrow$ is $f_1$ and $\searrow$ is $f_2$.}\label{Fig1}
\end{figure}

\begin{figure}[t]
\centerline{\includegraphics[width=16cm]{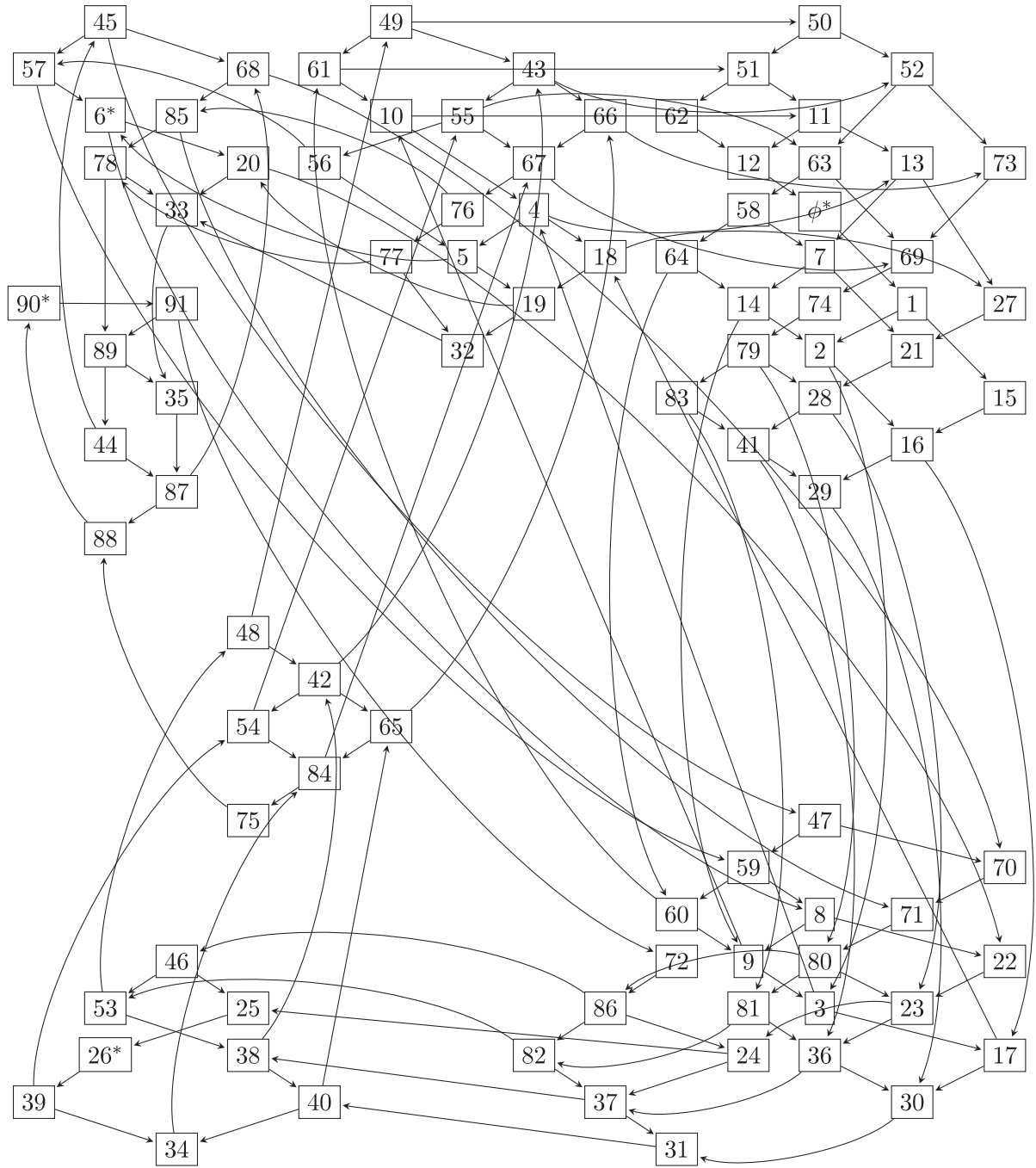}}

\caption{Crystal graph of $B_2$. $\searrow$ is $f_0$, $\swarrow$ is $f_1$ and others are $f_2$.}\label{Fig2}
\end{figure}

In this section we present crystal graphs of the $U'_q\big(G_2^{(1)}\big)$-crystals $B_1$ and $B_2$ in Figs.~\ref{Fig1}
and \ref{Fig2}. In the graphs $b\overset{i}{\longrightarrow}b'$ stands for $b'=f_ib$. Minimal elements are marked as~$*$.
Recall that as a~$U_q(G_2)$-crystal
\[
B_1\simeq B(0)\oplus B(\La_1),\qquad B_2\simeq B(0)\oplus B(\La_1)\oplus B(2\La_1).
\]
We give the table that relates the numbers in the crystal graphs to our representation of
elements according to which $U_q(G_2)$-components they belong to.

\medskip\noindent
$B(0):\;\fbox{$\phi^\ast$}=(0,0,0,0,0,0)$

\medskip\noindent
$B(\La_1)$:
\begin{alignat*}{4}
\fbox{1}&=(1,0,0,0,0,0)&
\fbox{2}&=(0,1,0,0,0,0)&
\fbox{3}&=\big(0, \tfrac{2}{3},\tfrac{2}{3},0,0,0\big)&
\fbox{4}&=\big(0, \tfrac{1}{3},\tfrac{4}{3},0,0,0\big)\\
\fbox{5}&=\big(0, \tfrac{1}{3},\tfrac{1}{3},1,0,0\big)&
\fbox{6$^\ast$}&=\big(0, \tfrac{1}{3},\tfrac{1}{3},\tfrac{1}{3},\tfrac{1}{3},0\big)&
\fbox{7}&=(0, 0,1,1,0,0)&
\fbox{8}&=\big(0, 0, 1, \tfrac{1}{3},\tfrac{1}{3},0\big)\\
\fbox{9}&=\big(0,0,0, \tfrac{4}{3},\tfrac{1}{3},0\big)&
\fbox{10}&=\big(0,0,0 , \tfrac{2}{3},\tfrac{2}{3},0\big)&
\fbox{11}&=(0,0,0,0,1,0)&
\fbox{12}&=(0,0,0,0,0,1)\\
\fbox{13}&=(0,0,2,0,0,0)&
\fbox{14}&=(0,0,0,2,0,0)
\end{alignat*}
$B(2\La_1)$:
\begin{alignat*}{4}
\fbox{15}&=(2,0,0,0,0,0)&
\fbox{16}&=(1,1,0,0,0,0)&
\fbox{17}&=\big(1, \tfrac{2}{3},\tfrac{2}{3},0,0,0\big)&
\fbox{18}&=\big(1, \tfrac{1}{3},\tfrac{4}{3},0,0,0\big)\\
\fbox{19}&=\big(1, \tfrac{1}{3},\tfrac{1}{3},1,0,0\big)&
\fbox{20}&=\big(1, \tfrac{1}{3},\tfrac{1}{3},\tfrac{1}{3},\tfrac{1}{3},0\big)&
\fbox{21}&=(1,0,1,1,0,0)&
\fbox{22}&=\big(1,0,1,\tfrac{1}{3},\tfrac{1}{3},0\big)\\
\fbox{23}&=\big(1,0,0,\tfrac{4}{3},\tfrac{1}{3},0\big)&
\fbox{24}&=\big(1, 0,0,\tfrac{2}{3},\tfrac{2}{3},0\big)&
\fbox{25}&=(1,0,0,0,1,0)&
\fbox{26$^\ast$}&=(1,0,0,0,0,1)\\
\fbox{27}&=(1,0,2,0,0,0)&
\fbox{28}&=(1, 0,0,2,0,0)&
\fbox{29}&=(0,2,0,0,0,0)&
\fbox{30}&=\big(0, \tfrac{5}{3},\tfrac{2}{3},0,0,0\big)\\
\fbox{31}&=\big(0, \tfrac{4}{3},\tfrac{4}{3},0,0,0\big)&
\fbox{32}&=\big(0, \tfrac{4}{3},\tfrac{1}{3},1,0,0\big)&
\fbox{33}&=\big(0, \tfrac{4}{3},\tfrac{1}{3},\tfrac{1}{3},\tfrac{1}{3},0\big)&
\fbox{34}&=(0,1,1,1,0,0)\\
\fbox{35}&=\big(0,1,1,\tfrac{1}{3},\tfrac{1}{3},0\big)&
\fbox{36}&=\big(0,1,0,\tfrac{4}{3},\tfrac{1}{3},0\big)&
\fbox{37}&=\big(0,1,0, \tfrac{2}{3},\tfrac{2}{3},0\big)&
\fbox{38}&=(0,1,0 ,0,1,0)\\
\fbox{39}&=(0,1,0,0,0,1)&
\fbox{40}&=(0,1,2,0,0,0)&
\fbox{41}&=(0,1,0,2,0,0)&
\fbox{42}&=\big(0, \tfrac{2}{3},\tfrac{2}{3},0,1,0\big)\\
\fbox{43}&=\big(0, \tfrac{1}{3},\tfrac{4}{3},0,1,0\big)&
\fbox{44}&=\big(0, \tfrac{1}{3},\tfrac{1}{3},1,1,0\big)&
\fbox{45}&=\big(0, \tfrac{1}{3},\tfrac{1}{3},\tfrac{1}{3},\tfrac{4}{3},0\big)&
\fbox{46}&=(0,0,1,1,1,0)\\
\fbox{47}&=\big(0,0,1, \tfrac{1}{3},\tfrac{4}{3},0\big)&
\fbox{48}&=\big(0,0,0, \tfrac{4}{3},\tfrac{4}{3},0\big)&
\fbox{49}&=\big(0,0,0,\tfrac{2}{3},\tfrac{5}{3},0\big)&
\fbox{50}&=(0,0,0,0,2,0)\\
\fbox{51}&=(0,0,0,0,1,1)&
\fbox{52}&=(0,0,2,0,1,0)&
\fbox{53}&=(0,0,0,2,1,0)&
\fbox{54}&=\big(0, \tfrac{2}{3},\tfrac{2}{3},0,0,1\big)\\
\fbox{55}&=\big(0, \tfrac{1}{3},\tfrac{4}{3},0,0,1\big)&
\fbox{56}&=\big(0, \tfrac{1}{3},\tfrac{1}{3},1,0,1\big)&
\fbox{57}&=\big(0, \tfrac{1}{3},\tfrac{1}{3},\tfrac{1}{3},\tfrac{1}{3},1\big)&
\fbox{58}&=(0, 0,1,1,0,1)\\
\fbox{59}&=\big(0,0,1,\tfrac{1}{3},\tfrac{1}{3},1\big)&
\fbox{60}&=\big(0,0,0,\tfrac{4}{3},\tfrac{1}{3},1\big)&
\fbox{61}&=\big(0,0,0,\tfrac{2}{3},\tfrac{2}{3},1\big)&
\fbox{62}&=(0,0,0,0,0,2)\\
\fbox{63}&=(0,0,2,0,0,1)&
\fbox{64}&=(0,0,0,2,0,1)&
\fbox{65}&=\big(0, \tfrac{2}{3},\tfrac{8}{3},0,0,0\big)&
\fbox{66}&=\big(0, \tfrac{1}{3},\tfrac{10}{3},0,0,0\big)\\
\fbox{67}&=\big(0, \tfrac{1}{3},\tfrac{7}{3},1,0,0\big)&
\fbox{68}&=\big(0, \tfrac{1}{3},\tfrac{7}{3},\tfrac{1}{3},\tfrac{1}{3},0\big)&
\fbox{69}&=(0,0,3,1,0,0)&
\fbox{70}&=\big(0,0,3,\tfrac{1}{3},\tfrac{1}{3},0\big)\\
\fbox{71}&=\big(0,0,2, \tfrac{4}{3},\tfrac{1}{3},0\big)&
\fbox{72}&=\big(0,0,2, \tfrac{2}{3},\tfrac{2}{3},0\big)&
\fbox{73}&=(0,0,4,0,0,0)&
\fbox{74}&=(0,0,2,2,0,0)\\
\fbox{75}&=\big(0, \tfrac{2}{3},\tfrac{2}{3},2,0,0\big)&
\fbox{76}&=\big(0, \tfrac{1}{3},\tfrac{4}{3},2,0,0\big)&
\fbox{77}&=\big(0, \tfrac{1}{3},\tfrac{1}{3},3,0,0\big)&
\fbox{78}&=\big(0, \tfrac{1}{3},\tfrac{1}{3},\tfrac{7}{3},\tfrac{1}{3},0\big)\\
\fbox{79}&=(0,0,1,3,0,0)&
\fbox{80}&=\big(0,0,1,\tfrac{7}{3},\tfrac{1}{3},0\big)&
\fbox{81}&=\big(0,0,0,\tfrac{10}{3},\tfrac{1}{3},0\big)&
\fbox{82}&=\big(0,0,0,\tfrac{8}{3},\tfrac{2}{3},0\big)\\
\fbox{83}&=(0,0,0,4,0,0)&
\fbox{84}&=\big(0, \tfrac{2}{3},\tfrac{5}{3},1,0,0\big)&
\fbox{85}&=\big(0, \tfrac{1}{3},\tfrac{4}{3},\tfrac{4}{3},\tfrac{1}{3},0\big)&
\fbox{86}&=\big(0,0,1,\tfrac{5}{3},\tfrac{2}{3},0\big)\\
\fbox{87}&=\big(0, \tfrac{2}{3},\tfrac{5}{3},\tfrac{1}{3},\tfrac{1}{3},0\big)&
\fbox{88}&=\big(0, \tfrac{2}{3},\tfrac{2}{3},\tfrac{4}{3},\tfrac{1}{3},0\big)&
\fbox{89}&=\big(0, \tfrac{1}{3},\tfrac{1}{3},\tfrac{5}{3},\tfrac{2}{3},0\big)&
\fbox{90$^\ast$}&=\big(0, \tfrac{2}{3},\tfrac{2}{3},\tfrac{2}{3},\tfrac{2}{3},0\big)\\
\fbox{91}&=\big(0, \tfrac{1}{3},\tfrac{4}{3},\tfrac{2}{3},\tfrac{2}{3},0\big)
\end{alignat*}

Comparing our crystal graphs with those in \cite{Y} we found that some 2-arrows are missing
in Fig.~3 of \cite{Y}.

\subsection*{Acknowledgements}

KCM thanks the faculty and staf\/f of Osaka University for their hospitality during his visit in August,
2009 and acknowledges partial support from NSA grant H98230-08-1-0080.
MM would like to thank Universiti Tun Hussein Onn Malaysia for supporting this study.
MO would like to thank the organizers of the conference ``Geometric Aspects of Discrete and
Ultra-Discrete Integrable Systems" held during March 30 -- April 3, 2009 at Glasgow for a warm
hospitality and acknowledges partial support from JSPS grant No. 20540016.

\newpage

\pdfbookmark[1]{References}{ref}
\LastPageEnding

\end{document}